\documentclass[12pt,leqno]{amsart}
\usepackage{amssymb,verbatim,enumerate,ifthen,amsmath}
\usepackage{mathtools}
\usepackage[mathscr]{eucal}
\usepackage[utf8]{inputenc}
\usepackage{hyperref}
\usepackage[T1]{fontenc}
\def\N{\mathbb{N}}
\def\B{\mathbb{B}}
\def\R{\mathbb{R}}
\def\V{\mathbb{V}}

\newtheorem{theorem}{Theorem}[section]
\newtheorem*{theorem*}{Theorem}
\def\Thm#1#2{\ifthenelse{\equal{#1}{*}}{\begin{theorem*}#2\end{theorem*}}
            {\begin{theorem}\label{T#1}#2\end{theorem}}}
\newtheorem{Atheorem}{Theorem}

\def\thm#1{Theorem~\ref{T#1}}
\newtheorem{proposition}[theorem]{Proposition}
\newtheorem*{proposition*}{Proposition}
\def\Prp#1#2{\ifthenelse{\equal{#1}{*}}{\begin{proposition*}#2\end{proposition*}}
{\begin{proposition}\label{P#1}#2\end{proposition}}}
\def\prp#1{Proposition~\ref{P#1}}
\newtheorem{corollary}[theorem]{Corollary}
\newtheorem*{corollary*}{Corollary}
\def\Cor#1#2{\ifthenelse{\equal{#1}{*}}{\begin{corollary*}#2\end{corollary*}}
{\begin{corollary}\label{C#1}#2\end{corollary}}}
\def\cor#1{Corollary~\ref{C#1}}

\newtheorem{lemma}[theorem]{Lemma}
\newtheorem*{lemma*}{Lemma}
\def\Lem#1#2{\ifthenelse{\equal{#1}{*}}{\begin{lemma*}#2\end{lemma*}}
             {\begin{lemma}\label{L#1}#2\end{lemma}}}

\theoremstyle{definition}
\newtheorem{remark}[theorem]{Remark}
\newtheorem*{remark*}{Remark}
\def\Rem#1#2{\ifthenelse{\equal{#1}{*}}{\begin{remark}\rm #2\end{remark}}
             {\begin{remark}\label{R#1}\rm #2\end{remark}}}

\newtheorem{example}[theorem]{Example}
\newtheorem*{example*}{Example}
\def\Exa#1#2{\ifthenelse{\equal{#1}{*}}{\begin{example*}\rm #2\end{example*}}
             {\begin{example}\label{Ex#1}\rm #2\end{example}}}

\def\eq#1{{\rm(\ref{E#1})}}
\def\Eq#1#2{\ifthenelse{\equal{#1}{*}}
  {\begin{equation*}\begin{aligned}#2\end{aligned}\end{equation*}}
  {\begin{equation}\begin{aligned}\label{E#1}#2\end{aligned}\end{equation}}}

\begin{document}
\begin{flushright}
\end{flushright}
\vspace{5mm}

\date{\today}

\title{Generalizing the concept of Bounded Variation}

\author[A. R. Goswami]{Angshuman R. Goswami}
\address{Department of Mathematics, University of Pannonia, 
Veszprém, Egyetem u. 10, H-8200.}
\email{goswami.angshuman.robin@mik.uni-pannon.hu}


\begin{abstract}
Let $[a,b]\subseteq\R$ be a non empty and non singleton closed interval and 
$P=\{a=x_0<\cdots<x_n=b\}$ is a partition of it. Then $f:I\to\R$ is said to be a function of $r$-bounded variation, if the expression $\overset{n}{\underset{i=1}{\sum}}|f(x_i)-f(x_{i-1})|^{r}$ is bounded for all possible partitions like $P$. One of the main result of the paper deals with the generalization of Classical Jordan decomposition theorem. We have shown that for $r\in]0,1]$, a function of $r$-bounded variation can be written as the difference of two monotone functions. While for $r>1$, under minimal assumptions such functions can be treated as approximately monotone function which can be closely approximated by a non decreasing majorant. We also proved that for $0<r_1<r_2$; the function class of $r_1$-bounded variation is contained in the class of functions satisfying $r_2$-bounded variations.

We go through approximately monotone functions and present a possible decomposition for $f:I(\subseteq \R_+)\to\R$ satisfying the functional inequality
$$f(x)\leq f(x)+(y-x)^{p}\quad (x,y\in I\mbox{  with  $x<y$  and  $ p\in]0,1[ $}).$$
A generalized structural study has also be done in that specific section.
 
On the other hand for $\ell[a,b]\geq d$; a function satisfying the following monotonic condition under the given assumption  will be termed as $d$-periodically increasing
$$f(x)\leq f(y)\quad \mbox{for all}\quad x,y\in I\quad\mbox{with}\quad y-x\geq d.$$

we establish that in a compact interval any bounded function can be decomposed as the difference of a monotone and a $d$-periodically increasing function.

The core details related to  past results, motivation, structure of each and every sections are thoroughly discussed below.
\end{abstract}
\maketitle
\section{Introduction}
The aim of this paper is to study the generalize notion of bounded variation. The paper can be subdivided into two parts. In the first part, we focused on power oriented bounded variation. While on the second half we studied distance dependent variation. At the beginning of each section, a brief discussion about structural characterization and various inclusion properties are also provided.\\
Let $r\in]0,\infty[$ be fixed. A function $f:[a,b]\to\R$ is said to be satisfying $r$-bounded variation if for any partition $P=\{a=x_0,\cdots ,x_n=b\}$ with $x_0<\cdots <x_n;$ the expression $\overset{n}{\underset{i=1}{\sum}}|f(x_i)-f(x_{i-1})|^{r}$ is bounded. It can be easily observed that for the case $r=1$, this notion implies the well-known definition of bounded variation. 

The term functions of bounded variation was first introduced by Jordan in his paper \cite{Jordan}. Later the concept was generalized by Wiener in the paper \cite{Wiener} where he investigated several topics related to Fourier Analysis. His generalized form is widely known $p$-variation and the definition of it is same as above mentioned $r$-variation with $r\geq 1$. Since then many more generalized forms in this direction were introduced including well known Young's generalized bounded variation in \cite{Young}.

One of our main objectives in this paper is to find a possible decomposition for the functions satisfying $r$-bounded variation.  We established that for $r\in ]0,1[$, if a function equipped with $r$-bounded variation then $f$ is nothing but a function of ordinary bounded variation. In other words, in such cases $f$ can be expressed as the difference of two monotone functions. We generalized the Jordan's decomposition theorem and  showed that a function $f:I(\subseteq R_+)\to\R$ with $0\in I$ satisfying Wiener variation/p-variation ($r$-variation with $r\geq 1$) can be treated as approximately monotone. Moreover with a minimal assumption on the error term, such a function $f$ can be precisely approximated by a nondecreasing function $g$ that satisfies the following inequality
$$\Phi(x)\leq g(x)-f(x)\leq 2\Phi\Big(\dfrac{x}{2}\Big)\quad\mbox{where}$$
\Eq{*}{
\Phi(u):=\sup_{x,x+u\in I}\Big\{ {^{r}}v_{f}(x+u)-{^{r}}v_{f}(x)\Big\}^{\dfrac{1}{r}}, \quad \quad \mbox{(${^{r}v}_{f}$ is $r$-variation function of $f$)} .
}
We also proved that for $0<r_1<r_2$, a function satisfying $r_1$-bounded variation is also a function of $r_2$-bounded variation. We showed that if $\Phi^{r}$ is superadditive, then in that scenario a $\Phi$-H\"older function satisfies $r$-bounded variation. The definition and other essentials for $\Phi$-H\"older functions are discussed below.

Let $p\geq 0$ be fixed and $f:I\to\R$ satisfies the following functional inequality
$$f(x)\leq f(y)+(y-x)^{p}\quad \mbox{for all $x,y\in I$} \quad \mbox{with $x<y$}.$$
For $p=1$, it is clearly visible that $f(t)+t$ is nondecreasing in $I$.  While for $p>1$, it can be shown that the function $f$ possesses usual monotonicity(increasingness). In the case $p=0$, $f$ can be treated as $\epsilon$-monotone function. Using the result related to the decomposition of such functions which was mentioned by Páles in his paper \cite{Monotonicity}; $f$ is expressible as $f:=g+h$, where $g$ is nondecreasing and $||h||\leq \dfrac{1}{2}$. The decomposition of $f$ when $p\in]0,1[$ is still unknown. In this paper we have shown that if $I\subseteq\R_+$ with $0\in I$, then for any $p\in]0,1[$, the function $f$ can be tightly approximated with a nondecreasing majorant $g$ that satisfies
$$f(x)+x^p\leq g(x)\leq f(x)+2^{(1-p)}x^p.$$

Motivated from this concept, we studied $\Phi$-monotonicity  in depth in our papers  \cite{Gos_Pal20} and \cite{Goswami}.
A function $f:I\to\R$ is said to be $\Phi$-monotone if for any $x,y\in I$ with $x<y$, the following  inequality holds
\Eq{*}{
f(x)\leq f(y)+\Phi(y-x)\quad\mbox{where}\quad\mbox{$\Phi:[0,\ell(I)]\to\R_+$ }
.}

On the other hand a function $f:I\to\R$ is termed as \emph{$\Phi$-H\"older} if both $f$ and $-f$ are $\Phi$-monotone, That is for any $x,y\in I$, the following inequality must satisfies 
\Eq{*}{
  |f(x)-f(y)|\leq\Phi(|y-x|).
}

Let $d>0$ be fixed and $\ell(I)\geq d$ such that for any $x,y\in I$ with $y-x \geq d$; $\Phi(y-x)=0$ holds. Then the notion of $\Phi$-monotonicity for $f$ can be treated as distance dependent monotone(increasing) function.

Inspired by this, we introduce the following concept. For a fixed $d>0$, a function $f: I\to \R$ is said to be increasing by a period $d$ (or $d$- periodically increasing) if the following  holds
$$f(x)\leq f(y) \qquad \mbox{ for all    }\quad x,y\in I\quad\mbox{with} \quad  y-x\geq d.$$
In our recently submitted paper, we have gone through the fundamental structural properties of such functions. For any given function $f$ we give a precise formula to obtain the largest $d$-periodically increasing minorant. We also show approximation of $d$-periodically increasing function by an ordinary monotone function. Besides basic characterization and a conditional decomposition of $d$-monotone function was also provided.

The main objective in the last section of the paper to establish a relationship between bounded functions and $d$ periodically increasing functions. We showed that a bounded function can be decomposed as the difference of a monotone and a $d$-periodically increasing function. For this establishment, we need to introduce a weaker concept of boundedness, named as $d$-variation.

Let $\ell([a,b])\geq d$ and the interval is partitioned as follows
$$P=\{a=x_0,\cdots,x_n=b\} \mbox{ such that $x_0<\cdots<x_n$ with $x_i-x_{i-1}\geq d$, $(i\in\{1,\cdots,n\}$)}.$$
Then for any function $f:[a,b]\to\R$ we define $d$-variation on the partition $P$ as 
\Eq{*}{\overset{b}{\underset{a}{V_{d}}}(f,P):=\sum_{i=1}^{n}|f(x_i-f(x_{i-1})|\quad \mbox{with}\quad x_i-x_{i-1}\geq d \quad (i\in\{1,\cdots,n\}).}
If for such possible partitions of $I$, the above expression remain bounded, we say that $f$ is a function of $d$-bounded variation. As mentioned before, we will later see that a function of $d$-bounded variation satisfies boundedness and vice verse. We also showed that the variation function $_{d}v_{f}$ is non decreasing. Moreover, superadditive characteristic of $d$-variation function is also shown.
\section{On approximately monotone and H\"older Functions}
Before proceeding, we need to recall several notions and terminologies. We consider $I$ be a non empty and non singleton interval. $\ell(I)(\in\overline{\R_+})$ denotes the length of the interval in the extended real line.
Then the function $\Phi:[0,\ell(I)]\to\R_+$ is said to be a error function. 

Based upon $\Phi$, we can formulate the concept of $\Phi$-monotonicity. A function $f: I\to\R$ is said to be $\Phi$-monotone(increasing) if for any $x,y\in I$ with $x\leq y$ the following functional inequality holds
$$f(x)\leq f(y)+\Phi(y-x).$$
On the other hand if both $f$ and $-f$ are $\Phi$-monotone, we will say that $f$ is $\Phi$-H\"older. In other words a function $f:I\to\R$ is said to be $\Phi$ H\"older if for any $x,y\in I$; it satisfies the inequality below
$$|f(x)- f(y)|\leq\Phi(|y-x|).$$
Now we will go through some inequalities associated with the error function $\Phi$ provided it has some structural properties .
\Prp{88}{Suppose $\Phi:[0,\ell(I)[\to\R_+$ be  subadditive and concave. Then for any $x,y\in[0,\ell(I)]$ with $x\leq y$, the following inequality holds 
\Eq{902}{\Phi(y)-\Phi(x)\leq\Phi(y-x)\leq 2\Phi\Big(\dfrac{y}{2}\Big)-\Phi(x).}
 }
Additionally, under the pre-mentioned conditions $(-\Phi)$ is $\Phi$-monotone.
\begin{proof}By using the subadditivity of $\Phi$, we get
\Eq{*}{
\Phi(y)=\Phi(y-x+x)\leq \Phi(y-x)+\Phi(x).
}
This is the first inequality of \eq{902}.

Now applying concavity property of $\Phi$, we can compute the following inequality
\Eq{*}{
\dfrac{\Phi(y-x)}{2}+\dfrac{\Phi(x)}{2}\leq\Phi\Big(\dfrac{(y-x)+x}{2}\Big)=\Phi\Big(\dfrac{y}{2}\Big).
}
Multiplying the above inequality by $2$ and rewriting it, we will obtain the remaining part of inequality of \eq{902}.

The second assertion can be easily seen from the first part of the inequality \eq{902}. By rearranging the terms of it we get
$$(-\Phi)(x)\leq (-\Phi)(y)+\Phi(y-x).$$
This yields the $\Phi$-monotonicity of $(-\Phi)$ and completes the proof of first assertion.
\end{proof}
Based upon the above theorem, we can propose the following corollary related to power functions.
\Cor{903}{For any $0\leq x\leq y$; the following inequality will hold true for any $p\in[0,1]$
\Eq{980}{y^p-x^p\leq(y-x)^{p}\leq 2^{1-p}y-x^p.}}
\begin{proof}
For $u\geq 0$ and for any fixed $p\in[0,1]$, the function $u\to u^{p}$ is subadditive and concave. Hence, the result directly follows from \prp{88} . 
\end{proof}
We are now ready to present a possible decomposition for $\Phi$-monotone function. The following theorem shows the tight approximation of a $\Phi$-monotone function by an ordinary monotone function. 
\Thm{777}{Let $I\subseteq\R_+$ be a non-empty and non-singleton interval such that $0\in I$. \\
$\Phi:[0,\ell(I)]\to\R_+$ is a subadditive error function. If $g:I\to\R$ is non-decreasing then $f:=g-\Phi$ is a $\Phi$-monotone function in $I$.
Conversely if, $f:I\to\R$ is $\Phi$-monotone where $\Phi$ is subadditive and concave. Then there exists a non decreasing majorant $g:I\to\R$ of $f$, that satisfies the following inequality
\Eq{1122}{\Phi(x)\leq g(x)-f(x)\leq 2\Phi\Big(\dfrac{x}{2}\Big).}}
\begin{proof}To prove the first part of the theorem, we assume $x,y\in I$ with $x\leq y$. By using the first inequality of \eq{902}, we get the following for $f$
\Eq{*}{
f(x)&=g(x)-\Phi(x)\\
&\leq g(y)-\Phi(y)+\big(\Phi(y)-\Phi(x)\big)\\
&=f(y)+\Phi(y-x).
}
This establishes that $f$ is $\Phi$-monotone.
To prove the converse assertion,
we use the inequality associated with $\Phi$-monotonicity of $f$. Applying the last inequality of \eq{902} there; we observe that the following holds
\Eq{555}{
f(x)+\Phi(x)\leq f(y)+2\Phi\bigg(\dfrac{y}{2}\bigg)\qquad ( \mbox{ for all } x,y\in I \,\,\mbox{ with }\, x\leq y).}
We can define $g:I\to\R$ as as follows
$$g(u):=\inf_{z\in I, u\leq z}\bigg(f(z)+2\Phi\Big(\dfrac{z}{2}\Big)\bigg).$$
We assume $x,y\in I$ with $x\leq y$. Then we can formulate the following inequality
\Eq{*}{
g(x)=\inf_{z\in I, x\leq z}\bigg(f(z)+2\Phi\Big(\dfrac{z}{2}\Big)\bigg)\leq \inf_{z\in I, y\leq z}\bigg(f(z)+2\Phi\Big(\dfrac{z}{2}\Big)\bigg)=g(y).
}
This shows that $g$ is nondecreasing. Now from \eq{555}, we obtain the following
\Eq{*}{
f(x)+\Phi(x)\leq g(x)\leq f(x)+2\Phi\bigg(\dfrac{x}{2}\bigg).
}
Subtracting $f$ from all part of above inequality, we will get \eq{1122} and it validates the result.
\end{proof}

From the above theorem, we can establish a conditional decomposition of the function $f: I\to\R$ satisfying the following inequality 
\Eq{bass}{f(x)\leq f(y)+c(y-x)^{p} \quad x,y\in I\quad\mbox{with} \quad x\leq y. \quad (p\in]0,1[).}
We heve already discussed about such type of functions in the Introduction Section of this paper.
\Cor{2233}{Let $I\subseteq\R_+$ be a nonempty and nonsingleton interval with $0\in I$. Then for $f:I\to\R$ satisfying \eq{bass}, there exists a nondecreasing majorant $g$ of $f$ such that the following inequality holds in $I$.
$$cx^p\leq g(x)-f(x)\leq \big(2^{1-p}c\big)x^p.$$}
\begin{proof}
One can easily show that for a fixed $c>0$, the function $u\to c u^{p}$ ($u\in I$) is subadditive and concave. Hence, the result is a direct implication of \thm{777} and \cor{903}.
\end{proof}
In the next section, we will discuss about  functions of $r$-bounded variation.

\section{On functions of power dependent variation}
Through out this section $I$ will denote the non empty and non singleton interval $[a,b]$. Let $P$ be a partition of $I$ as follows
\Eq{001}{
P=\{x_0,\cdots,x_n\} \quad \mbox{with} \quad a=x_0 <\cdots <x_n=b}
Let $r\in]0,\infty[$ be fixed. For a real valued function  $f: I\to\R$, we define $r$-bounded variation $\overset{b}{\underset{a}{V}^{r}}(f,P)$ for the partition $P$ of $I$ as below
$$\overset{b}{\underset{a}{V}^{r}}(f,P):=\sum_{i=0}^{n}|f(x_i)-f(x_{i-1})|^{r}.$$ 
If the above mentioned expression is bounded with respect to all such possible partitions of $I$, we say that $f$ possesses $r$-bounded variation. The supremum of all $r$-variations of $f$ is called total $r$-variation which is denoted by  
 $$\overset{b}{\underset{a}{V^{r}}}(f):=\sup \overset{b}{\underset{a}{V}^{r}}(f,P).$$
We term ${^{r}v}_{f}:I\to\R$, defined by ${^{r}v}_{f}(x):=\overset{x}{\underset{a}{V^{r}}}(f)$ as $r$-variation function for $f$.
It can be clearly observed that if a function is unbounded then it can't be a function of $r$-bounded variation. Also, there are bounded functions which can not be categorized as function satisfying $r$-bounded variation for any $r\in]0,\infty[$. For example 
 $f:I\to\R$ defined by
$$f(x)=\begin{cases}
1, \quad\mbox{if $x$ is rational}\\
0 \quad\mbox{if $x$ is irrational}
\end{cases}
$$
is bounded but not a function of $r$-bounded variation.
We will utilize the following subadditive (superadditive) property in the upcoming result extensively.

For any $r\leq 1$ and $a_1,\cdots,a_n\in\R_+$, the following inequality satisfies
\Eq{*}{
 \bigg(\sum_{i=1}^{n} a_i\bigg)^r\leq\sum_{i=1}^{n} a_i^{r}.
}
For $r\geq 1$ under the same assumptions on ${a_{i}}'s$, the reserve inequality holds.

Now we are going to show that the $r$- variation function possesses monotonicity.
\Prp{91}{${^{r}v}_{f}$ is a nondecreasing function.}
\begin{proof}Let $x,y\in I$ with $x\leq y$. For any $\epsilon>0$, we will have a partition $P$ of $I$ as
defined in \eq{001} satisfying the following inequality 
\Eq{*}{{^{r}v}_{f}(x)&< \sum_{i=0}^{n}|f(x_i)-f(x_{i-1})|^{r}+\epsilon\\
& \leq \sum_{i=0}^{n}|f(x_i)-f(x_{i-1})|^{r}+|f(y)-f(x)|^{r}+\epsilon\\
&\leq {^{r}v}_{f}(y)+\epsilon.}
Upon taking $\epsilon\to 0$, we get the monotonicity of ${^{r}v}_{f}$.
\end{proof}

The next proposition will show that total $r$-bounded variation of a function defined in a compact interval of $\R$ possesses characteristics similar to superadditivity.
\Prp{9991}{Let $r>0$ be fixed and $f:I\to\R$ is a function of $r$- bounded variation. Then for any $c\in I^{\circ}$, the following inequality  holds 
\Eq{1000}{\overset{c}{\underset{a}{V}^{r}}(f)+\overset{b}{\underset{c}{V}^{r}}(f)\leq \overset{b}{\underset{a}{V}^{r}}(f).}
Additionally, in case of $r\in]0,1];$ the above inequality turns into equality.}
\begin{proof}
To establish the inequality, we assume that $\epsilon>0$ is arbitrary. Then there must exists two partitions $P_1=\{a=x_0,\cdots,x_k=c\}$ and $P_2=\{c=x_k,\cdots,x_n=b\}$ respectively for $[a,c]$ and $[c,b]$ with $x_0<\cdots <x_n$ satisfying the following two inequalities
\Eq{*}{
\overset{c}{\underset{a}{V}^{r}}(f)<\overset{c}{\underset{a}{V}^{r}}(f,P_1)+\dfrac{\epsilon}{2}\quad\mbox{and}\quad \overset{c}{\underset{a}{V}^{r}}(f)<\overset{b}{\underset{c}{V}}(f,P_2)+\dfrac{\epsilon}{2}.}
Upon adding up the two inequalities side by side, we obtain
\Eq{*}{\overset{c}{\underset{a}{V}^{r}}(f)+\overset{b}{\underset{c}{V}^{r}}(f)&<\overset{c}{\underset{a}{V}^{r}}(f,P_1)+\overset{c}{\underset{a}{V}^{r}}(f,P_2)+\epsilon\\
&=\overset{b}{\underset{a}{V}^{r}}(f,P_1\cup P_2)+\epsilon\\
&\leq \overset{b}{\underset{a}{V}}(f)+\epsilon.
}
Since $\epsilon$ is arbitrary, taking $\epsilon\to 0$, we obtain \eq{1000}. This validates establishment the first statement.\\
To prove the second assertion it will be enough to show that for any $r\in]0,1]$, the reverse inequality of \eq{1000} also holds. For any arbitrary $\epsilon>0$, there must exists a partition like $P$ of $I$, as defined in \eq{001} satisfying the following inequality
\Eq{*}{
\overset{b}{\underset{a}{V}^{r}}(f)<\sum_{i=1}^{n}|f(x_i)-f(x_{i-1})|^r+\epsilon.
}
Then for any $c\in I^{\circ}$ there must exists a $k\in\{1,\cdots, n\}$ such that $c\in [x_{k-1},x_k]$ where $x_{k-1},x_k\in P$. Together with this, by using the subadditive property for $r\in]0,1]$; we can expand the above inequality as follows
\Eq{*}{\overset{b}{\underset{a}{V}^{r}}(f)&<\sum_{i=1}^{k-1}|f(x_i)-f(x_{i-1})|^r+|f(x_{k-1}-f(c)+f(c)-f(x_k)|^{r}\\
&\qquad \qquad \qquad \qquad \qquad \qquad \qquad+\sum_{i=k-1}^{n}|f(x_i)-f(x_{i-1})|^r+\epsilon\\
&<\sum_{i=1}^{k-1}|f(x_i)-f(x_{i-1})|^r+\Big(|f(x_{k-1}-f(c)|+|f(c)-f(x_k)|\Big)^{r}\\
&\qquad \qquad \qquad \qquad \qquad \qquad \qquad+\sum_{i=k-1}^{n}|f(x_i)-f(x_{i-1})|^r+\epsilon\\
&\leq \sum_{i=1}^{k-1}|f(x_i)-f(x_{i-1})|^r+|f(x_{k-1}-f(c)|^{r}+|(f(c)-f(x_{k+1})|^{r}\\
&\qquad \qquad \qquad \qquad \qquad \qquad \qquad+\sum_{i=k}^{n}|f(x_i)-f(x_{i-1})|^r+\epsilon\\
&\leq \overset{c}{\underset{a}{V}^{r}}(f)+\overset{b}{\underset{c}{V}^{r}}(f)+\epsilon.}
Since $\epsilon$ is arbitrary, considering it negligibly small we obtain 
\Eq{*}{\overset{b}{\underset{a}{V}^{r}}(f)\leq  \overset{c}{\underset{a}{V}^{r}}(f)+\overset{b}{\underset{c}{V}^{r}}(f).}
This proves the validity of assertion and completes the proof.
\end{proof}
The next proposition will show that that algebraic multiplication of two $r$-bounded function results a function with same variation characteristic provided $r\in]0,1].$
\Prp{inc}{Let $r_1,r_2\in ]0,1]$. If $f,g:I\to\R$ are two functions of $r_1$ and $r_2$-bounded variations  respectively. Then $fg:I\to\R$ will satisfy $\max\{r_1,r_2\}$-bounded variation. }
\begin{proof}
We assume $r=\max\{r_1,r_2\}$. Since both $f$ and $g$ are functions of $r_1$ and $r_2$-bounded variations respectively; then by \thm{last}, $f$ and $g$ will satisfy $r$-bounded variation.  Therefore there exists a $k\in\R_+$ such that both the inequalities $|f(x)|\leq k$ and $|g(x)|\leq k$ holds for all $x\in I$. Now, we consider a partition $P$ as in \eq{001}. Since $r\in]0,1]$, utilizing subadditive property of $V^r\big(P,(fg)\big)$, we compute upper bound as follows
\Eq{*}{
V^r\big(P,(fg)\big)&=\sum_{i=1}^n|(fg)(x_i)-(fg)(x_{i-1})|^r\\
&=\sum_{i=1}^n|f(x_i)g(x_i)-f(x_i)g(x_{i-1})+f(x_i)g(x_{i-1})-f(x_{i-1})g(x_{i-1})|^r\\
&\leq\sum_{i=1}^n\Big(|f(x_i)||g(x_i)-g(x_{i-1})|+|f(x_i)-f(x_{i-1})||g(x_{i-1})|\Big)^r\\
&\leq k^r\bigg(\sum_{i=1}^{n}|f(x_i)-f(x_{i-1})|^r+\sum_{i=1}^{n}|g(x_i)-g(x_{i-1})|^r\bigg)\\
&\leq k^r\Big[V^r\big(P,f\big)+V^r\big(P,g\big)\Big]\\
&\leq k^r\Big[V^r\big(f\big)+V^r\big(g\big)\Big].
}
We will ended up the same bound as above for any arbitrary partition $P$. And hence it is evident that $fg$ satisfies $r$-bounded variation. This completes the proof of the statement.
\end{proof}
However, the above statement is not valid for $r>1$.
We consider $I=[0,1]$ and take $r>1$ be fixed. Now we construct two functions $g,h:I\to\R$ as follows
\Eq{*}{
g(x):=
x\cos\Big(\dfrac{\pi}{x}\Big)  
\quad \quad \mbox{and} \quad \quad h(x):=\begin{cases}0 \quad \mbox{if $x=0$}\\
\dfrac{1}{x} \quad \mbox{ if $x\in]0,1]$}
\end{cases}}
One can easily check that both $g$ and $h$ are functions of $r$-bounded variation. Now the function $f:=g.h$  can be defined as follows 
\Eq{*}{
f(x):=\begin{cases}
0 \quad \mbox{if x=0}\\
\cos\Big(\dfrac{\pi}{x}\Big) \quad \mbox{if $x\in]0,1]$}
\end{cases}.}
It is evident that $f$ is not a function of $r$-bounded variation for any $r>0$. Thus we can conclude that depending upon the value of $r$, the structural composition of $r$-bounded variation varies.

Now we will show the relationship between a function satisfying $r$-bounded variation and approximately H\"older function. But before that we need to go through the proposition below
\Prp{phi}{ For $r\in[1,\infty[$; the function $\Phi:[0,\ell(I)]\to\R_+$ is defined as follows
\Eq{phi}{
\Phi(u):=\sup_{x,x+u\in I}\Big\{ {^{r}}v_{f}(x+u)-{^{r}}v_{f}(x)\Big\}^{\dfrac{1}{r}}.
}
Then $\Phi$ is subadditive function.
}
\begin{proof}
To proof the statement we assume $u,v\in\R_+$ such that $u+v\in \ell(I).$ Then from \eq{phi}, for any given $\epsilon> 0$ there must exists a $x\in I$ with $x+u+v\in I$ satisfying the following functional inequality which is extended using subadditive property due to $\dfrac{1}{r}.$
\Eq{*}{
\Phi(u+v)&< \Big\{{^{r}}v_{f}(x+u+v)-{^{r}}v_{f}(x)\Big\}^{\dfrac{1}{r}}+\epsilon\\
&=\Big\{{^{r}}v_{f}(x+u+v)-{^{r}}v_{f}(x+u)\Big\}^{\dfrac{1}{r}}+\Big\{{^{r}}v_{f}(x+u)+{^{r}}v_{f}(x)\Big\}^{\dfrac{1}{r}}+\epsilon\\
&\leq \sup_{x,x+v\in I}\Big\{ {^{r}}v_{f}(x+v)-{^{r}}v_{f}(x)\Big\}^{\dfrac{1}{r}}+\sup_{x,x+u\in I}\Big\{ {^{r}}v_{f}(x+u)-{^{r}}v_{f}(x)\Big\}^{\dfrac{1}{r}}+\epsilon\\
&=\Phi(u)+\Phi(v)+\epsilon.
}
Since $\epsilon$ is arbitrary, upon taking $\epsilon\to 0$, it is evident from the above inequality that $\Phi$ possesses subadditivity.
\end{proof}
\Prp{000}{
Let $r\in[1,\infty[$ and the function $f:I\to\R$ satisfies $r$-bounded variation. Then $f$ is a $\Phi$-H\"older function where $\Phi$ is defined as in \eq{phi}.
}
\begin{proof}
To prove the theorem we assume $x,y\in I$ with $x<y$. Assuming $u=y-x$ and utilizing \prp{9991} and \prp{phi} we compute the following inequality
\Eq{*}{
|f(x)-f(y)|\leq \Bigg(\overset{y}{\underset{x}{V}^{r}}(f)\Bigg)^{\dfrac{1}{r}}\leq \Big\{{^{r}}v_{f}(y)-{^{r}}v_{f}(x)\Big\}^{\dfrac{1}{r}}\leq \Phi(y-x).
}
This validates our statement.
\end{proof}
The next result is a generalization of Jordan's decomposition theorem. We will see that depending upon the value of $r$; a $r$-bounded function can be expressed as the difference of two monotone functions or under some minimal assumptions can be approximated by a nondecreasing function. 
\Thm{233}{Let $r\in]0,\infty[$ be fixed. $f: I\to\R$ be a function of $r$-bounded variation. 
If $r\in]0,1]$, $f$ can be written as the difference of two monotone functions. Additionally, for $r=1$; the inverse implication also holds true.\\
On the other hand, for $r\in]1,\infty[$; if $\inf(I)=0$ and $\Phi$ is concave then $f$ can be approximated by a monotone majorant $g$ satisfying the same functional inequality as in \eq{1122}. That is
$$\Phi(x)\leq g(x)-f(x)\leq 2\Phi\bigg(\dfrac{x}{2}\bigg),$$
where the function $\Phi$ is as defined in \eq{phi}.}
\begin{proof}
First we will consider the case for $r\in]0,1]$.\\
We can express $f$ as $f:=\big({^{r}v}_{f}\big)^{\dfrac{1}{r}}-\Big(\big({^{r}v}_{f}\big)^{\dfrac{1}{r}}-f\Big)$. Let, $x,y\in I$ with $x<y$. Due to non-negativity of ${^{r}v}_{f}$ and by \prp{91}; it is evident that the function $\big({^{r}v}_{f}\big)^{\dfrac{1}{r}}$ is monotone. Now to prove the first  assertion, we only to show that $\big({^{r}v}_{f}\big)^{\dfrac{1}{r}}-f$ is also nondecreasing. Let $x,y\in I$ with $x<y$. Then from \prp{9991}, we obtain the following inequality
\Eq{*}{
|f(y)-f(x)|^r\leq \overset{y}{\underset{x}{V}^{r}}(f)={^{r}v}_{f}(y)-{^{r}v}_{f}(x).
}
Since $r\in]0,1]$, by using superadditivity and non negativity of ${^{r}v}_{f}$, above inequality can be expanded as follows
\Eq{*}{
f(y)-f(x)\leq |f(y)-f(x)|\leq \Big({^{r}v}_{f}(y)-{^{r}v}_{f}(x)\Big)^{\dfrac{1}{r}}\leq \big({^{r}v}_{f}(y)\big)^{\dfrac{1}{r}}-\big({^{r}v}_{f}(x)\big)^{\dfrac{1}{r}}.
} 
Rearranging the terms in the above expression we obtain
\Eq{*}{
\big({^{r}v}_{f}(x)\big)^{\dfrac{1}{r}}-f(x)\leq \big({^{r}v}_{f}(y)\big)^{\dfrac{1}{r}}-f(y).}
This yields the monotonicity of $\big({^{r}v}_{f}\big)^{\dfrac{1}{r}}-f$ and completes the first establishment.

For r=1, the function $f\in\B\V(I)$ and hence by Jordan's Decomposition Theorem the inverse implication is obvious.\\
Finally, we arbitrarily select $r\in]1,\infty]$. From \prp{000}, we have that $f$ is a $\Phi$-H\"older function. Hence $f$ satisfies $\Phi$-monotonicity. Again from \prp{phi} and our assumption, we have that $\Phi$ is a concave and  subadditive error function. Thus by the converse part of \thm{777} the inequality holds. And this completes the theorem. 
\end{proof}

However the reverse implication of the above theorem does not hold.  
Based on above results, we are now going to study some inclusion properties of different variation classes. The class of functions satisfying $r$-bounded variation where $r\in]0,1[$ will be denoted by $\underset{]0,1[}{\B\V}(I)$. The notation ${\B\V}(I)$ is used for the class of functions of ordinary bounded variation. On the other hand, for $r>1$, such class is represented by $\underset{]1,\infty[}{\B\V}(I)$. Based upon the definition of these classes we can establish the following inclusion theorem. One can also opt the below result as a direct consequence of \thm{last}. From \prp{inc}and the theorem below, one can easily see that the class of functions $\underset{]0,1[}{\B\V}(I)$ is a ring under usual functional addition and multiplication. While $\underset{]1,\infty[}{\B\V}(I)$ is a commutative group.  
\Thm{BV}{ Let $r\in]0,1[$, then the following proper inclusion holds
\Eq{09}{\underset{]0,1[}{\B\V}(I)\subsetneq{\B\V}(I)\subsetneq\underset{]1,\infty[}{\B\V}(I).}}
\begin{proof}
To establish the first inclusion property, we assume that $f \in \underset{]0,1[}{\B\V}(I).$
Then from \thm{233}, we know that $f$ can be decomposed as the algebraic difference of two monotone functions. Utilizing Jordan's decomposition theorem; we can conclude $f \in {\B\V}(I)$. This shows the first part of inclusion in \eq{09}.\\
Now to show ${\B\V}(I)\subseteq\underset{]1,\infty[}{\B\V}(I)$, we assume $f\in{\B\V}(I)$ and $r\in]1,\infty[$ arbitrary. Then at the partition $P$ as defined in \eq{001}; utilizing superadditivity, we calculate $V^r(f,P)$ as follows
\Eq{*}{
V^r(f,P)=\sum_{i=1}^{n}|f(x_i)-f(x_{i-1}|^r\leq \bigg(\sum_{i=1}^{n}|f(x_i)-f(x_{i-1})|\bigg)^r
.}
Since the expression $\overset{n}{\underset{i=1}{\sum}}|f(x_i)-f(x_{i-1}|$ is bounded, this implies boundedness of $V^r(f,P)$ as well. Likewise  same holds true for any other partition of $I$.
This implies $f\in \underset{]1,\infty[}{\B\V}(I)$ and completes the inclusion property.

Now, we are going to show that the inclusions are proper through examples. First we show that exists a function $f\in {\B\V}$ such that $f\notin \underset{]0,1[}{\B\V}$. Similarly we also construct example such that $f\in \underset{]1,\infty[}{\B\V}(I)\setminus {\B\V}(I)$.\\
The identity function ($t\to t$) is a function of bounded variation in $I.$ We assume $r\in]0,1[$ be fixed and consider  a partition of $I$ as mentioned in \eq{001} with the condition  
 $x_i-x_{i-1}=\frac{b-a}{n}$ for all $i\in\{1,\cdots n\}.$  we calculate  $V^{r}\big(f,P\big)$ as follows
\Eq{*}{
V^{r}\big(f,{P}\big)&=\sum_{i=1}^{n}|x_i-x_{i-1}|^{r} =n^{(1-r)}(b-a)^{r}.}
Upon taking $n\to\infty$, we can see that $V^{r}\big(f,{P}\big)$ is not bounded. Through arbitrariness of $r$, we conclude that the identity function does not satisfy $r$-bounded variation for any $r<1$. It shows the validity of our first assertion.\\
To justify the second assertion we consider $I$ as $[0,1]$ and define the following function $f:[0,1]\to\R$ as follows
\Eq{*}{
f(x)=\begin{cases}
x\cos\Big(\dfrac{\pi}{x}\Big) \mbox{ if }x\in]0,1]\\
0 \quad \mbox{ if } \quad  x=0.
\end{cases}
}
We choose, $m\in\N$  such that $\frac{1}{m}<1.$ Then for any $n\in N$; We partition $[0,1]$ as \\
$\hat{P}=\Big\{0,\dfrac{1}{m+n},\dfrac{1}{m+n-1},\cdots\dfrac{1}{m},1\Big\}.$ 
 Based on the partition $\hat{P}$, we calculate the ordinary variation of $f$ as follows
\Eq{*}{
V^{1}(f,\hat{P})\geq\Big(\dfrac{1}{m+n}+\dfrac{1}{m+n-1}\Big)+\cdots+\Big(\dfrac{1}{m}+\dfrac{1}{m+1}\Big)
}
Upon taking $n\to\infty$, it is evident that $V^{1}(f,\hat{P})$ is not bounded. So, we can conclude $f\notin {\B\V}(I).$\\
Now we assume $r>1$ be arbitrary. Then we compute $V^{r}(f,\hat{P})$ to obtain the following inequality
\Eq{*}{
V^{r}(f,\hat{P})&\leq \Big(1+\dfrac{1}{m+n}\Big)^{r}+\Big(\dfrac{1}{m+n}+\dfrac{1}{m+n-1}\Big)^{r}+\cdots+\Big(\dfrac{1}{m}+\dfrac{1}{m+1}\Big)^{r}+\Big(1+\dfrac{1}{m}\Big)^{r}\\
&\leq \Big(1+\dfrac{1}{m+n}\Big)^{r}+2^{r}\Big(\dfrac{1}{m+n-1}\Big)^{r}\cdots+2^{r}\Big(\dfrac{1}{m}\Big)^{r}+\Big(1+\dfrac{1}{m}\Big)^{r}
}
Since for any $r>1$; the series $\underset{n\in \N}{\sum}\bigg(\dfrac{1}{n}\bigg)^{r}$ is convergent, then by comparison theorem, we have that $V^{r}(f,\hat{P})$ is bounded. The same assertion can be shown for any other partition of $[0,1]$ just by approximating that partition with $\hat{P}$. This implies $f$ is a function of $r$-bounded variation for any $r>1$. In other words $f\in\underset{]1,\infty[}{\B\V}.$ This validates that the second inclusion is proper and completes the proof.
\end{proof}
The next proposition shows the linkage between $\Phi$-H\"older function and $r$-bounded variation. 
\Prp{99}{Let $\Phi: [0,\ell(I)]\to\R_+$ be an error function such that ${\Phi}^{r}$  is superadditive. Then any $\Phi$-H\"older function will also be a function of $r$-bounded variation.} 
\begin{proof}
 $f:I\to\R$ be a $\Phi$ -H\"older function. we partition $P$ as described in \eq{001}. Now by using the superadditivity of  ${\Phi}^{r}$, we obtain the $r$-variation of $f$ as follows
\Eq{*}{
\underset{a}{\overset{b}{V^{r}}}(f,P)&=\sum_{i=1}^{r}|f(x_i)-f(x_{i-1})|^{r}\\
&\leq\sum_{i=1}^{n}\Phi^{r}(x_i-x_{i-1})\\
&\leq \Phi^{r}(b-a).}
This holds true with respect to any possible partition of $I$ and hence $f$ satisfies $r$- bounded variation. 
\end{proof}
For the next proposition we assume that $\B\V_{r}$  represents classes of functions satisfying $r$- bounded variation. The inclusion statement of \thm{BV} is a direct implication of the following theorem. 
\Thm{last}{ Let $0<r_1<r_2$. Then
$\B\V_{r_1}\subseteq \B\V_{r_2}$ holds.
}
\begin{proof}
To prove the statement we assume that 
$f\in\B\V_{r_1}.$ We consider a partition $P$ as defined in \eq{001}. We compute $V^{r_2}\big(f,P\big)$ as follows.
\Eq{*}{V^{r_2}\big(f,P\big):=\sum_{i=1}^{n}|f(x_i)-f(x_{i-1})|^{r_2}}
Since $\dfrac{r_2}{r_1}>1$, using the superadditivity of non negative numbers we compute the following inequality
\Eq{*}{V^{r_2}\big(f,P\big)=\sum_{i=1}^{n}\bigg(|f(x_i)-f(x_{i-1})|^{r_1}\bigg)^{\frac{r_2}{r_1}}\leq \bigg(\sum_{i=1}^{n}|f(x_i)-f(x_{i-1})|^{r_1}\bigg)^{\frac{r_2}{r_1}}\leq\bigg(V^{r_1}\big(f)\bigg)^{\frac{r_2}{r_1}}.}
This shows that the above inequality also implies boundedness of  $V^{r_2}\big(f,P\big)$. Since $P$ is any arbitrary partition, we can conclude $f\in\B\V_{r_2}$. This completes the proof.
\end{proof}
Through example we can easily observe that the above inclusion is proper. For the justification here we only consider the case where $0<r_1<r_2<1$. One can easily notice the same validation for the remaining cases. For simplicity we consider $I=[0,1]$ and select any $r>1$. Now we define the function $f:\R\to\R$ as follows
\Eq{examp}{
f(x)=\begin{cases}
(-1)^n\dfrac{1}{n^{r}}, \mbox{ if }x=\dfrac{1}{n} \mbox{ for all }n\in\N\\
0 \quad  \mbox{ otherwise }
\end{cases}
}
Since $r>1$ the series $\underset{n\in \N}\sum\bigg(\dfrac{1}{n^r}\bigg)$ is convergent. Using it, one can easily show that $f\in{\B\V}$. Now we select $r_1,r_2\in ]0,1[$ such that $\dfrac{1}{r_2}<r<\dfrac{1}{r_1}$. Now we choose a fixed $m\in\N$ and consider the partition $\hat{P}$ for any arbitrary $n\in \N$ as mentioned in \thm{BV}
$$\hat{P}=\bigg\{0,\dfrac{1}{m+n},\dfrac{1}{m+n-1}\cdots \dfrac{1}{m},1\bigg\}.$$
By calculating $r_1$-variation of $f:I\to\R$ at $\tilde{P}$ we obtain the following inequality
\Eq{*}{
V^{r_1}(f,\hat{P})
&\geq 2^{rr_1}\Big(\dfrac{1}{m+n-1}\Big)^{rr_1}+\cdots+2^{rr_1}\Big(\dfrac{1}{m+1}\Big)^{rr_1}
}
As $rr_1<1$; the series $\underset{n\in \N}{\sum}\bigg(\dfrac{1}{n^{rr_1}}\bigg)$ is divergent. Thus upon taking $n\to\infty$ we can clearly notice that $V^{r_1}(f,P)$ is unbounded. Therefore $f$ is not a $r_1$-bounded function.
In a analogous way we can easily establish that $f\in \B\V_{r_2}$. This justifies the statement.

\section{On functions of $d$-bounded variation}
Throughout this section $d$ will be a fixed positive number and $I=[a,b]$ will represent a non-empty and non singleton closed interval with $\ell(I)\geq 2d$. However for some of the results this condition can be relaxed to $\ell(I)\geq d$ and it is very much self explanatory. A function $f: I\to \R$ is said to be  $d$-periodically increasing if for any $x,y\in I$, with $y-x\geq d$, $f(x)\leq f(y)$ holds. It is evident that the class of $d$-periodically increasing functions is a convex cone as it is closed under addition and multiplication with non negative scalars. Moreover if $f$ is $d$-periodically increasing and non negative then  for any $r(\in \R_+)$ the function $f^{r}$ is also periodically increasing with the same period $d$. It can also be shown that the class of $d$-periodically increasing functions  is also closed under pointwise limit operations.

Let $f: I\to\R$ be a real valued function. We consider an arbitrary partition of $I$ as
$P=\{x_0,\cdots,x_n\}$ with $ a=x_0 <\cdots <x_n=b$ satisfying $ x_i-x_{i-1}\geq d$  where $i\in\{1,\cdots,n\}$. 
Based on the partition $P$, we derive the $d$-variation of $f$ in $I$ as follows
\Eq{dvar}{\overset{b}{\underset{a}{V_{d}}}(f,P):=\sum_{i=0}^{n}|f(x_i)-f(x_{i-1})|.} 
The supremum of \eq{dvar} with respect to all such $d$-variation of $f$, is termed as total $d$-variation of $f$.  
In other words, we can define total $d$-variation of $f$ as follows 
\Eq{*}{\overset{b}{\underset{a}{V_{d}}}(f):=\sup\overset{b}{\underset{a}{V_{d}}}(f,P).} 
By using the notion of total $d$-variation of $f$ we can formulate the $d$-variation function $_{d}v_{f}$ as follows
\Eq{*}{_{d}v_{f}(x):=
\begin{cases}
0\quad \qquad \mbox{if}\quad x\in[a,a+d[\\
\overset{x}{\underset{a}{V_{d}}}(f)\quad\mbox{if}\quad x\in[a+d,b]
\end{cases}}
It is evident that any bounded function $f$ will satisfy $d$- bounded variation since
$${\overset{b}{\underset{a}{V_{d}}}}(f,P)\leq  \bigg\lceil{\dfrac{b-a}{d}}\bigg\rceil\sup_{u,v\in[a,b]} |f(u)-f(v)|$$
holds. Moreover from the inequalities
\Eq{*}{|f(x)|-|f(b)|\leq |f(b)-f(x)|\leq {\overset{b}{\underset{a}{V_{d}}}}(f) \quad \mbox{for all}\quad x\in[a,b-d] \\
|f(x)|-|f(a)|\leq |f(x)-f(a)|\leq {\overset{b}{\underset{a}{V_{d}}}}(f) \quad \mbox{for all}\quad x\in[a+d,b] 
}
and from our pre assumption on $I$ that $\ell(I)\geq 2d$, it can be easily seen that 
\Eq{258}{|f(x)|\leq \max\big\{ |f(a)|,|f(b)|\big\}+{\overset{b}{\underset{a}{V_{d}}}}(f).
}
Which yields that
any function $f: I\to\R$ satisfying $d$-bounded variation is bounded as well.  
Therefore, we can use the terms bounded function and function satisfying $d$-bounded variation interchangeably. 
The following results give a structural overview of this newly defined variation. The proposition below shows the $d$-periodic monotonicity of a $d$- variation function.
\Prp{223}{Let $f: I\to\R$ be a bounded function. Then $d$-variation function $_{d}v_{f}$ is $d$-periodically increasing.}
\begin{proof}
To validate the statement we consider two cases. first let $x\in[a,a+d[$ and $y\in[a+d,b]$ such that $y-x\geq d$. Then the $d$-periodical monotonicity of ${_{d}v}_{f}$ is clearly evident from the following inequality
$${_{d}v}_{f}(x)=0\leq \overset{y}{\underset{a}{V_{d}}}(f)={_{d}v}_{f}(y).$$
Next we assume $x,y\in[a+d,b]\subseteq I$ such that $y-x\geq d$. Let $\epsilon\geq 0$ be arbitrary. Then there must exists a partition $P$ of $[a,x]$ as $P=\{a=x_0,\cdots ,x_n=x\}$ such that $a=x_0<x_1<\cdots<x_n=x$ with $x_i-x_{i-1}\geq d$ for all $i\in\{1,\cdots n\}$ satisfying
$$ {_{d}v}_{f}(x)=\overset{x}{\underset{a}{V_{d}}}(f)< \overset{x}{\underset{a}{V_{d}}}(f,P)+\epsilon.$$
Using the above inequality we obtain the following
\Eq{*}{\overset{y}{\underset{a}{V_{d}}}(f)+\epsilon
&\geq\bigg(\sum_{i=1}^{n}|f(x_i)-f(x_{i-1})|+|f(y)-f(x)|\bigg)+\epsilon\\
&=\overset{x}{\underset{a}{V_{d}}}(f,P)+\epsilon+|f(y)-f(x)|\\
&> \overset{x}{\underset{a}{V_{d}}}(f).}
Taking $\epsilon\to 0$; we see that $\overset{x}{\underset{a}{V_{d}}}(f)\leq \overset{y}{\underset{a}{V_{d}}}(f)$ holds. 
It proves that $_{d}v_{f}$ is $d$-periodically increasing in the interval $I$.
\end{proof}

Next proposition shows that the total $d$-variation of a function possesses superadditive property in restricted interval.

\Prp{113}{Let $f$ be a function of $d$-bounded variation in $I$. Then, for any $c\in [a+d, b-d]\subseteq I$; the following inequality holds
$$\overset{c}{\underset{a}{V_{d}}}(f)+\overset{b}{\underset{c}{V_{d}}}(f)\leq \overset{b}{\underset{a}{V_{d}}}(f).$$}
\begin{proof}
Since at the beginning of the section we assumed that $\ell(I)\geq 2d$, therefore $[a+d,b-d]$ is non empty. Let $c\in [a+d,b-d]\subseteq I$; this implies that $\ell([a,c])$ and $\ell([c,b])\geq d$. We assume $\epsilon>0$ be arbitrary. Then there must exists two partitions for $[a,c]$ and $[c,b]$ respectively as $P_1=\{x_0,\cdots,x_{r-1},c=x_r\}$
and $P_2=\{c=x_r, x_{r+1},\cdots,x_n\}$ such that $a=x_0<\cdots<x_r(=c)<\cdots<x_n=b$ with $x_i-x_{i-1}\geq d$ for all $i\in\{1,\cdots,n\}$ satisfying the following two inequalities
\Eq{*}{
\overset{c}{\underset{a}{V_{d}}}(f)<\overset{c}{\underset{a}{V_{d}}}(f,P_1)+\dfrac{\epsilon}{2}\qquad\mbox{and}\qquad
\overset{b}{\underset{c}{V_{d}}}(f)<\overset{b}{\underset{c}{V_{d}}}(f,P_2)+\dfrac{\epsilon}{2}
}
Summing up these two inequalities side by side,we obtain
\Eq{*}{
\overset{c}{\underset{a}{V_{d}}}(f)+\overset{b}{\underset{c}{V_{d}}}(f)&<\overset{c}{\underset{a}{V_{d}}}(f,P_1)+\overset{b}{\underset{c}{V_{d}}}(f,P_2)+\epsilon\\
&=\overset{b}{\underset{a}{V_{d}}}(f,P_1\cup P_2)+\epsilon\\
&\leq \overset{b}{\underset{a}{V_{d}}}(f)+\epsilon.
}
Due to arbitrariness of $\epsilon$, upon taking $\epsilon\to 0$, we obtain
\Eq{}{
\overset{c}{\underset{a}{V_{d}}}(f)+\overset{b}{\underset{c}{V_{d}}}(f)\leq \overset{b}{\underset{a}{V_{d}}}(f)
.}
This yields the inequality to be shown and validates the statement.
\end{proof}
We will use the following proposition in the result related to decomposition. 
\Prp{11}{Every $d$-periodically increasing function satisfies $d$-bounded variation. }
\begin{proof}
 We assume $f:I\to\R$ is a $d$-periodically increasing function.
Then for the partition $P$ as defined in \eq{dvar}, the following equation can be compute
\Eq{*}{
\overset{b}{\underset{a}{V_{d}}}(f,P)=\sum_{i=1}^{n}|f(x_i)-f(x_{i-1})|=\sum_{i=1}^{n}\Big(f(x_i)-f(x_{i-1}\Big)=f(b)-f(a)
.}
Upon taking supremum of all such possible partitions $P$, we can conclude that
$$\overset{b}{\underset{a}{V_{d}}}(f)=f(b)-f(a).$$
It shows that $f$ is a function of $d$-bounded variation and completes the proof.
\end{proof}

Now we are ready to state a possible decomposition for any bounded function in terms of a monotone and $d$-periodically increasing functions. The result somewhat resembles with the classical Jordan decomposition theorem. 
\Thm{2222}{ A function of $d$-bounded variation can be expressed as the difference of a mononotone and a $d$-periodically increasing function.\\
Conversely difference of two $d$-periodically increasing functions is a function of $d$-bounded variation.}
\begin{proof}
First we assume that $f:I\to \R$ is a function of $d$-bounded variation. 
Now we construct a function ${_{d}w}_f(x):I\to\R_+$ as
$${_{d}w}_f:=
\begin{cases}
0\qquad \mbox{if}\quad x\in[a,a+d[\\
\underset{z\in [a, x-d]}{\sup}\overset{x}{\underset{z}{V_{d}}}(f)\,\, \mbox{if}\quad  x\in[a+d, b]
\end{cases}
$$
and rewrite $f$ in the following way
$$f:={_{d}}w_{f}-({_{d}}w_{f}-f).$$
To validate the theorem it will be enough to show that ${_{d}}w_{f}$ is increasing while ${_{d}}w_{f}-f$ is $d$ periodically increasing.\\
Let $x,y\in I$ with $x<y$. If $x,y\in[a,a+d[$ then there is nothing to show. For $x\in[a,a+d[$ and $y\in[a+d,b]$, the monotonicity follows from the definition of ${_{d}}w_{f}$. On the other hand if both $x,y\in[a+d,b]$, then the following inequality yields that ${_{d}}w_{f}$ is non decreasing
$${_{d}}w_{f}(x)=\underset{z\in [a, x-d]}{\sup}\overset{x}{\underset{z}{V_{d}}}(f)\leq \underset{z\in [a, y-d]}{\sup}\overset{y}{\underset{z}{V_{d}}}(f)={_{d}w}_f(y).$$ 

To show ${_{d}}w_{f}-f$ is $d$-periodically increasing, we consider $x,y\in I$ such that $y-x\geq d.$ This implies $y\in [a+d, b]$ and we have to consider consider two cases. For $x\in [a, a+d[$, by definition ${_{d}}w_{f}(x)=0$.
On the other hand if $x\in [a+d, b]$ then ${_{d}}w_{f}(x)=\underset{z\in [a, x-d]}{\sup}\overset{x}{\underset{z}{V_{d}}}(f)$. For both of these cases we can compute the following inequality
$${_{d}}w_{f}(y)\geq {_{d}}w_{f}(x)+|f(y)-f(x)|\geq {_{d}}w_{f}(x)+f(y)-f(x).$$
Rearranging the terms of the above inequality we see that ${_{d}}w_{f}-f$ is $d$-periodically increasing. This proves the first part of the theorem.

To show the converse assertion, let $g$ and $h$ are two $d$-periodically increasing functions and hence by \prp{11}, both of these functions satisfies $d$-bounded variation. Therefore, from \eq{258}, we have that both $g$ and $h$ are bounded. Thus $g-h$ will also be bounded.  In other words $g-h$  is a function of $d$-bounded variation. This completes the proof of the theorem.
\end{proof}
The occurrence of $d$-periodically increasing function often noticeable in the stock market prices of a reputed company's  share. In most of the cases it varies in a short time frame but eventually  it gets a upward trend after a distinct time interval. The similar type of periodical monotonicity is observable in the population growth graph of most East European nations. Therefore, interdisciplinary study of such delayed monotonicity is also very much possible.
\bibliographystyle{plain}

\end{document}